\documentstyle[amscd]{amsart}

\def\NN{{\blb N}}

\def\PP{{\blb P}}

\def\RR{{\blb R}}

\def\ZZ{{\blb Z}}

\def\blb#1{\Bbb#1}

\def\11{{1\kern-3.5pt 1}}
\def\mumu{{\mu\kern-4.2pt\mu}}

\def\boxtimes{\setbox0\hbox{$\Box$}\copy0\kern-\wd0\hbox{$\times$}}

\def\pol{{\bar p}}

\def\Mol{{\bar M}}

\def\Uol{{\overline{ U}}}

\def\End{\operatorname {End}}

\def\Hom{\operatorname {Hom}}
\def\im{\operatorname {im}}

\def\ker{\operatorname {ker}}

\def\Spec{\operatorname {Spec}}

\def\Ann{\operatorname{Ann}}

\def\dim{\operatorname{dim}}

\def\End{\operatorname{End}}

\def\Fdim{{\sf Fdim}}

\def\fl.dim{\operatorname{flat.dim}}
\def\Fract{\operatorname{Fract}}

\def\Gr{\operatorname{Gr}}
\def\GrMod{{\sf GrMod}}

\def\Hom{\operatorname{Hom}}

\def\id{\operatorname{id}}

\def\max{\operatorname{max}}

\def\Mod{{\sf Mod}}

\def\Proj{{\sf Proj}}

\def\rank{\operatorname{rank}}

\def\Spec{\operatorname{Spec}}

\def\Tails{{\sf Tails}}

\def\d{\downarrow}

\let\oldtext\text
\def\text#1{\oldtext{\normalshape #1}}

\def\a{\alpha}

\def\c{\gamma}
\def\d{\delta}

\def\ve{\varepsilon}

\def\s{\sigma}

\def\fp{{\frak p}}

\def\sT{{\sf T}}

\def\Qcoh{{\sf{ Qcoh}}}

\def\cE{{\cal E}}
\def\cF{{\cal F}}

\def\cI{{\cal I}}
\def\cK{{\cal K}}
\def\cL{{\cal L}}
\def\cM{{\cal M}}

\def\cO{{\cal O}}

\def \E {\operatorname{E}}

\newtheorem{lemma}{Lemma}[section]
\newtheorem{proposition}[lemma]{Proposition}
\newtheorem{theorem}[lemma]{Theorem}
\newtheorem{corollary}[lemma]{Corollary}

{

}

\theoremstyle{definition}

\newtheorem{example}[lemma]{Example}
\newtheorem{definition}[lemma]{\sl Definition}

\theoremstyle{remark}

\numberwithin{equation}{section}

\begin{document}

\pagenumbering{arabic}

\title{Integral non-commutative spaces}

\author{ S. Paul Smith}

\address{Department of Mathematics, Box 354350, Univ.
Washington, Seattle, WA 98195, USA}

\email{smith@@math.washington.edu}

\thanks{The author was supported by NSF grant DMS-0070560}


\begin{abstract}
A non-commutative space $X$ is a Grothendieck category $\Mod X$.
We say $X$ is integral if there is an indecomposable 
injective $X$-module $\cE_X$ such that its endomorphism ring 
is a division ring and every $X$-module is a subquotient of a 
direct sum of copies of  $\cE_X$.
A noetherian scheme is integral in this sense 
if and only if it is integral in the usual sense.
We show that several classes of non-commutative spaces are
integral.
We also define the function field and generic point of an integral
space and show that these notions behave as one might expect.
\end{abstract}

\maketitle

\section{Introduction}

We follow Rosenberg and Van den Bergh in taking a Grothendieck
category as our basic non-commutative geometric object. We think of
a Grothendieck category $\Mod X$ as ``the quasi-coherent sheaves on
an imaginary non-commutative space $X$''.  The commutative model
is the category $\Qcoh X$ of quasi-coherent sheaves on a 
quasi-separated, quasi-compact scheme $X$. 
The two non-commutative models are $\Mod R$, the category of right
modules over a ring, and $\Proj A$, the non-commutative projective
spaces defined by Verevkin \cite{Ver} and Artin and Zhang \cite{AZ}.

This paper defines $X$ to be {\sf integral} if $\Mod X$ is locally
noetherian and there is an indecomposable 
injective $X$-module $\cE_X$ such that $\End \cE_X$ is a division
ring and every $X$-module is a subquotient of a 
direct sum of copies of  $\cE_X$ (Definition \ref{defn.integral}). 
If $X$ is integral, then up to
isomorphism there is only one indecomposable injective with 
these properties.  The {\sf function field} of an
integral space is the division ring $\End \cE_X$.
We also define the generic point of an integral space. 
Corollary \ref{cor.integral.integral} shows that a
noetherian scheme is integral in the usual sense if and only if
$\Qcoh X$ is integral in our sense. In that case 
$\cE_X$ is the constant sheaf with sections equal to the function
field of $X$, and the function field in our sense coincides with 
the usual function field of $X$. 

Goldie's theorem implies that an affine
space having a prime right noetherian coordinate ring is integral.
However, we give a categorical definition of integrality 
so that it can be applied to those non-commutative spaces that
are not defined in terms of a ringed space. 
The non-commutative projective planes 
defined by Artin, Tate, and Van den Bergh
\cite{ATV2} are integral. The non-commutative analogues of 
$\PP^n$ associated to enveloping algebras of Lie algebras \cite{LV}, 
and the analogues of $\PP^n$ arising from the Sklyanin algebras 
\cite{TV} are integral. 
The exceptional fiber in Van den Bergh's blowup of a
non-commutative surface at a point \cite{vdB} is always integral.


Section five shows that non-commutative integral spaces 
enjoy some of the properties of integral schemes.

\medskip

{\bf Acknowledgements.}
The notion of integrality grew out of earlier work with J. Zhang. 
Theorem \ref{prop.ModM} is due to him, and
we thank him for allowing us to include it here.
We thank D. Happel, K. Goodearl, P. J\o rgensen, 
C. Pappacena, and J. Zhang for several helpful conversations.

We are grateful to the referee of an earlier version of this paper. He
suggested a change to an earlier definition of integrality, and 
that change represents a substantial improvement.

This work was begun while the author was visiting the
University of Copenhagen, and was continued during
the workshop on Noncommutative Algebra at the Mathematical
Sciences Research Institute at Berkeley.
The author is grateful to both organizations for
their hospitality and financial support.

\section{Preliminaries}
\label{sect.qsch}

Throughout we work over a fixed commutative base ring $k$. 
All categories are assumed to be $k$-linear, and so are all functors 
between them.

We adopt the framework for non-commutative algebraic geometry 
originated  by Rosenberg \cite{Rosen} and further developed by
Van den Bergh \cite{vdB}.
Definitions of terms we do not define can be found in \cite{vdB}.

\begin{definition}
A {\sf non-commutative space} $X$ is a Grothendieck category $\Mod X$.
Objects in $\Mod X$ are called $X$-modules.
We say $X$ is {\sf locally noetherian} if $\Mod X$ is locally 
noetherian (that is, if it has a set of noetherian generators).
\end{definition}

\begin{definition}
If $X$ and $Y$ are non-commutative spaces, 
a {\sf weak map} $f:Y \to X$
is a natural equivalence class of left exact 
functors $f_*:\Mod Y \to \Mod X$.
A weak map $f:Y \to X$ is a {\sf map} if $f_*$ has a left 
adjoint.
A left adjoint to $f_*$ will be denoted by $f^*$, and a right
adjoint will be denoted by $f^!$. 
\end{definition}

We say $X$ is {\sf affine} if $\Mod X$ has a progenerator, and in this case 
any ring $R$ for which $\Mod X$ is equivalent to $\Mod R$ is called a 
{\sf coordinate ring} of $X$.

If $(X,\cO_X)$ is a scheme then the category $\Mod \cO_X$ 
of all sheaves of $\cO_X$-modules is a Grothendieck category. 
If $X$ is quasi-compact and quasi-separated (for example, if $X$ is
a noetherian scheme) 
the full subcategory of $\Mod \cO_X$ consisting of the 
quasi-coherent $\cO_X$-modules is a Grothendieck category
\cite[page 186]{SGA6}.
We denote this category by $\Qcoh X$. Whenever $X$ is a
quasi-compact and quasi-separated scheme we will speak of it as a
space in our sense with the tacit understanding that 
$\Mod X$ is synonomous with $\Qcoh X$.

\section{Integral spaces, generic points, and function fields}
\label{sect.integral}

Throughout this section we fix a locally noetherian space $X$.

We denote the injective envelope of an $X$-module $M$ by $E(M)$.

\begin{definition}
\label{defn.integral}
A locally noetherian space $X$
is {\sf integral} if there is an indecomposable injective 
$\cE_X$ such that $\End \cE_X$ is a division ring and 
every $X$-module is a subquotient of a direct sum of copies of $\cE_X$. 
We call $\cE_X$ the {\sf big injective} in $\Mod X$.
\end{definition}

{\bf Remarks.}
The endomorphism ring of an indecomposable injective $\cE$
is a division ring if and only if $\Hom_X(\cE/N,\cE)=0$ for all
non-zero submodules $N$ of $\cE$.

When $X$ is locally noetherian the following conditions on an
$X$-module $\cE$ are equivalent: (a) every $X$-module is a 
subquotient of a direct sum of copies of $\cE$;
(b) every noetherian $X$-module is a subquotient of a finite 
direct sum of copies of $\cE$. 

Corollary \ref{cor.unique.inj} shows that the big
injective is unique up to isomorphism, thus justifying the use of
the definite article. Therefore the rank of a
module, the generic point, and the function field of $X$, all of
which are defined below in terms of $\cE_X$, are unambiguously defined.

\begin{definition}
Let $X$ be an integral locally noetherian space.
An $X$-module $M$ is {\sf torsion} if $\Hom(M,\cE_X)=0$.
A module is {\sf torsion-free} if the only
submodule of it that is torsion is the zero submodule.
\end{definition}


The torsion modules form a localizing subcategory of $\Mod X$.

\begin{definition}
Let $X$ be an integral locally noetherian space.
The {\sf rank} of an $X$-module $M$ is the length of
$\Hom_X(M,\cE_X)$ as a left $\End \cE_X$-module. 
We denote it by $\rank M$.
\end{definition}

Thus an $X$-module is torsion if and only if its rank is zero.

Because $\cE_X$ is injective, rank is additive on short exact
sequences.

The hypotheses on $\cE_X$ ensure that it has rank one, and
every proper quotient of it has rank zero. Hence every non-zero
submodule of $\cE_X$ has rank one.

Because a noetherian $X$-module is a subquotient of a finite direct
sum of copies of $\cE_X$, its rank is finite.

If $\rank M\ge 1$, then $M$ has a quotient of rank one, namely
$M/\ker f$ where $f$ is a non-zero element of $\Hom_X(M,\cE_X)$.

If $M$ is a noetherian torsion-free module of rank $n \ge 1$, 
then there is a finite chain $M=M_0 \supset M_1 \supset \ldots 
\supset M_{n-1} \supset M_{n}=0$ such that each 
$M_i/M_{i+1}$ is torsion-free of rank one. To see this begin by
choosing $M_1$ to be maximal subject to the condition that
$\rank(M_0/M_1)=1$; the maximality ensures that
$M_0/M_1$ is torsion-free, then argue by induction on $n$.

Since rank is additive on exact sequences, it induces
a group homomorphism $ \rank :K_0(X) \to \ZZ.$

\begin{lemma}
\label{lem.ModW}
Let $X$ be an integral locally noetherian space.
Let $M$ be a noetherian $X$-module. There exist 
noetherian submodules $L_1,\ldots,L_n$ of $\cE_X$,
a submodule $L \subset L_1 \oplus \ldots \oplus L_n$, and an
epimorphism $\varphi:L \to M$ such that $\varphi(L \cap L_i) \ne 0$
for all $i$.

Furthermore, the rank of $L$ is $n$.
\end{lemma}
\begin{pf}
By the definition of integrality there are 
noetherian submodules $L_1,\ldots,L_n$ of $\cE$,
a submodule $L \subset L_1 \oplus \ldots \oplus L_n$, and an
epimorphism $\varphi: L\to M$.
Choose this data so that $n$ is as small as possible. 
If $\varphi(L \cap L_i)$ were equal to zero, then there would
be an epimorphism $L/L \cap L_i \to M$, and since $L/L \cap
L_i$ is isomorphic to a submodule  of $L_1 \oplus \ldots \oplus
L_n/L_i$ this would contradict the minimality of $n$. So we conclude
that $\varphi(L \cap L_i) \ne 0$ for all $i$.

Since the rank of each $L_i$ is one, $\rank(L_1 \oplus \ldots \oplus
L_n) =n$. Thus $\rank L \le n$. However, $L \cap L_i \ne 0$ for
all $i$, whence $\rank L=n$.
\end{pf}

\begin{proposition}
\label{prop.red+irr}
Let $X$ be an integral locally noetherian space. 
If $J$ is a non-zero injective, then 
$\Hom_X(\cE_X,J) \ne 0$.
\end{proposition}
\begin{pf}
If $J$ is a non-zero injective $X$-module, then it contains a
non-zero noetherian submodule, say $N$.
Let $\varphi:L \to N$ be an epimorphism as in Lemma \ref{lem.ModW}.
The restriction of $\varphi$ to $L \cap L_1$, which is a submodule
of $\cE_X$, extends to a non-zero map map from $\cE_X$ to $J$.
\end{pf}

\begin{proposition}
\label{prop.torsion}
Let $X$ be an integral locally noetherian space.
An essential extension of a torsion module is torsion.
\end{proposition}
\begin{pf}
Let $P \subset M$ be an essential extension of a torsion module
$P$.
It suffices to prove the result when $M$ is noetherian because 
every $M$ is a directed union of noetherian
submodules $M_i$ each of which is an essential extension of $M_i \cap
P$.

Choose an epimorphism $\varphi:L \to M$  as in Lemma
\ref{lem.ModW}. Since $\varphi(L \cap L_i) \ne 0$, 
$P \cap \varphi(L \cap L_i) \ne 0$. But $P$ is torsion and
$L \cap L_i$ is torsion-free, so the restriction of $\varphi$ to $L
\cap L_i$ is not monic. Thus $\ker \varphi \cap L_i \ne 0$. Since
$L_i$ is torsion-free of rank one, $L_i/\ker \varphi \cap L_i$ is
torsion. Since $M$ is a subquotient of 
$\oplus_{i=1}^n L_i/\ker \varphi \cap L_i$ it is also a torsion
module.
\end{pf}

\begin{corollary}
\label{cor.unique.inj}
If $X$ is an integral locally noetherian space, there is only
one indecomposable injective up to isomorphism having the
properties in Definition \ref{defn.integral}.
\end{corollary}
\begin{pf}
Let $\cE_X$ be the injective in Definition \ref{defn.integral}, and
let $\cE$ be another indecomposable 
injective such that its endomorphism ring is a division
ring and every $X$-module is a subquotient of a direct sum of
copies of $\cE$. 

By Proposition \ref{prop.red+irr},
there is a non-zero map $\a:\cE_X \to \cE$. If $\a$ is monic, 
then its image would be a direct summand of $\cE$, so would equal 
$\cE$ because $\cE$ is
indecomposable; hence the result is true if $\a$ is monic.
Suppose to the contrary that $\a$ is not monic. 
Then its image is a proper quotient of $\cE_X$ so is
torsion. Therefore $\cE$ is the injective envelope of a torsion
module, so is itself torsion by Proposition \ref{prop.torsion}.
That is, $\Hom_X(\cE,\cE_X)=0$. It follows that $\Hom_X(-,\cE_X)$
vanishes on all $X$-modules. 
This is absurd, so we conclude that $\a$
is monic.
\end{pf}

\begin{definition}
\label{defn.fn.field}
Let $X$ be an integral locally noetherian space.
The {\sf function field} of $X$ is the division algebra 
$$
k(X):=\Hom_X(\cE_X,\cE_X).
$$
The {\sf generic point} of $X$ is the space
$\eta$ defined by
$$
\Mod \eta = \Mod X/\sT,
$$
where $\sT$ is the full subcategory consisting of the torsion
modules.
\end{definition}

Since $\sT$ is a localizing subcategory of $\Mod X$,
there is an adjoint pair of functors $(j^*,j_*)$ where 
$j^*:\Mod X \to \Mod \eta:=\Mod X/\sT$ is
the quotient functor, and $j_*$ its right adjoint. 
This defines a map of spaces 
$$
j:\eta \to X.
$$
For the rest of this section $j$ will denote this map.

\begin{proposition}
\label{prop.gen.pt}
Let $X$ be an integral locally noetherian space. 
If $\eta$ is its generic
point, then $\Mod \eta$ is equivalent to $\Mod \, k(X)$.
\end{proposition}
\begin{pf}
Since $\cE_X$ is torsion-free and every proper quotient of it is
torsion, $j^*\cE \cong j^*M$  for every non-zero submodule $M$
of $\cE_X$. It follows that $j^*\cE_X$ is a simple module 
in $\Mod \eta$.

If $M$ is an $X$-module, then  $E(M)/M$ is torsion by Proposition
\ref{prop.essl.extn}, so $j^*M \cong j^*E(M)$. Since $E(M)$ is a
direct sum of indecomposable injectives, and $j^*$ commutes with
direct sums, and an indecomposable injective is either torsion or
isomorphic to $\cE_X$, $j^*M$ is
isomorphic to a direct sum of copies of $j^*\cE_X$.
Therefore every $\eta$-module is isomorphic to 
a direct sum of copies of $j^*\cE_X$. 
Thus $\Mod \eta$ is equivalent to $\Mod D$ where
$D$ is the endomorphism ring of $j^*\cE_X$.

Since $\cE_X$ is torsion-free and injective, 
$j_*j^*\cE_X \cong \cE_X$, whence
$$
D=\Hom_{\eta}(j^*\cE_X,j^*\cE_X) \cong \Hom_X(\cE_X,j_*j^*\cE_X) 
\cong k(X).
$$
This completes the proof.
\end{pf}

{\bf Remark.}
The rank of an $X$-module $M$ 
is equal to the length of $j^*M$ as a right $k(X)$-module.
To see this, first observe that this length is equal to the length of
the left $k(X)$-module
$\Hom_{\eta}(j^*M,j^*\cE_X)$; second, observe that
we have the following natural isomorphisms:
\begin{align*}
\Hom_{\eta}(j^*M,j^*\cE_X) 
& \cong \Hom_{\eta}(j^*j_*j^*M,j^*\cE_X)
\\
& \cong \Hom_{X}(j_*j^*M,j_*j^*\cE_X)
\\
& \cong \Hom_{X}(j_*j^*M,\cE_X).
\end{align*}
It follows that the length of $j^*M$ is equal to the rank of
$j_*j^*M$.
However, there is an exact sequence $0 \to A \to M \to j_*j^*M \to
B \to 0$ where $A$ and $B$ are torsion modules, so $\rank M=\rank
j_*j^*M.$

\begin{theorem}
[Zhang]
\label{prop.ModM}
\label{thm.ModW}
Let $X$ be an integral locally noetherian space.
Then 
\begin{enumerate}
\item{}
every torsion-free module has a non-zero submodule that is 
isomorphic to a submodule of $\cE_X$;
\item{}
a uniform torsion-free module has rank one;
\item{}
the injective envelope of every torsion-free module of rank one is
isomorphic to $\cE_X$;
\item{}
$\cE_X$ is the unique indecomposable injective of rank one;
\item{}
every simple $X$-module is a subquotient of $\cE_X$.
\end{enumerate}
\end{theorem}
\begin{pf}
(1)
It suffices to prove this for a noetherian torsion-free module $M$.
Choose an epimorphism $\varphi:L \to M$ as in Lemma \ref{lem.ModW}.
Since $M$ is torsion-free and every proper quotient of 
$L \cap L_1$ is torsion, $\varphi(L
\cap L_1) \cong L \cap L_1$ which is a non-zero submodule of $\cE_X$.

(2)
It suffices to prove this for a noetherian torsion-free 
uniform module $M$.
Choose an epimorphism $\varphi:L \to M$ as in Lemma \ref{lem.ModW}
and set $M_i=\varphi(L \cap L_i)$. Thus $M_i$ is torsion-free of
rank one.
Since $M$ is uniform, $\cap_{i=1}^n M_i \ne 0$. 
An induction argument shows that the rank of $M_1+\ldots
+M_n$ is one: certainly $\rank(M_j)=1$ for all $j$, and 
\begin{align*}
\rank(M_1+\ldots +M_{i+1})
& =\rank(M_1+\ldots +M_{i})+ \rank(M_{i+1})
\\
& \qquad \qquad - \rank(M_1+\ldots +M_{i})\cap M_{i+1})
\\
& =\rank(M_1+\ldots +M_{i})+1-1
\\
&
=\rank(M_1+\ldots +M_{i}).
\end{align*}
But the rank of $L/\oplus_{i=1}^n (L \cap L_i)$ is zero, so
the rank of its quotient $M/\sum_{i=1}^n M_i$ is also zero.
Hence $\rank M=1$.

(3)
By (1) a rank one torsion-free module is an essential extension of a
non-zero submodule of $\cE_X$, so embeds in $\cE_X$.

(4)
Let $E'$ be an indecomposable injective of rank one.
Since $\E'$ is the injective envelope of all its non-zero submodules,
it follows from Proposition \ref{prop.torsion} that $E'$ is
torsion-free. Hence by (1) $E'$ and $\cE_X$ have a common
submodule, whence $E' \cong \cE_X$.

(5)
Let $S$ be a simple module. There is an epimorphism 
$\varphi:L \to S$ as in Lemma \ref{lem.ModW}.
There is a finite descending chain 
$L=K_0 \supset K_1 \supset \ldots \supset K_r=0$ of submodules 
such that each $K_i/K_{i+1}$ is torsion-free of rank one. 
Since $\Hom_X(L,S) \ne 0$, $\Hom_X(K_i/K_{i+1},S) \ne 0$
for some $i$. Since $S$ is simple, this provides the 
required epimorphism.
\end{pf}

The next result improves on Proposition \ref{prop.torsion}.

\begin{proposition}
\label{prop.essl.extn}
Let $X$ be a locally noetherian integral space.
If $L \subset M$ is an essential extension of $X$-modules, then 
$\rank L = \rank M$.
\end{proposition}
\begin{pf}
It is sufficient to prove the result when $M$ is the injective
envelope of $L$.
In that case, we can write $M$ as a direct sum of indecomposable
injectives, say $M =\oplus_i M_i$. Then $L \cap M_i \ne 0$ for all
$i$, and $M/L$ is a quotient of $ \oplus M_i/L \cap M_i$. Since
$M_i$ is an indecomposable injective, either its rank is zero or it
is isomorphic to $\cE_X$; in either case, $\rank M_i/L \cap M_i=0$.
Hence $\rank M/L=0$, and the result follows.
\end{pf}



\section{Examples of integral spaces}
\label{sect.egs}

A scheme $X$ is integral in the usual sense of 
algebraic geometry if and only if $\cO_X(U)$ is an integral domain
for all open subsets $U \subset X$. Corollary
\ref{cor.integral.integral} shows that
a noetherian scheme is integral in our sense  if and
only if it is integral in the usual sense.

We show that an affine space
having a prime right noetherian coordinate ring is integral.
We give other examples 
which indicate that our notion of integral is reasonable. In
particular, Theorem \ref{thm.intgl.tails} implies that 
the non-commutative analogues of $\PP^2$ discovered by
Artin-Tate-Van den Bergh are integral spaces, as are the Sklyanin
analogues of $\PP^n$.

\begin{proposition}
\label{prop.intgl.scheme}
Let $X$ be an integral noetherian scheme. 
Let $\cK$ denote the constant sheaf having sections the function
field of $X$.
If $\cM$ is a coherent $\cO_X$-module, then there is a 
coherent $\cO_X$-submodule, $\cL$ say, of a finite direct sum of copies
of $\cK$ and an epimorphism $\psi: \cL \to \cM$.
\end{proposition}
\begin{pf}
For the purposes of this proof we call a coherent $\cO_X$-module
$\cM$ good if there is such an epimorphism.
Clearly a finite direct sum of good modules is good, a submodule of
a good module is good, and a quotient of a good module is good.

Let $E(\cM)$ denote the injective envelope in $\Qcoh X$ of an
$\cO_X$-module. This is a direct sum of indecomposable injectives.
Each indecomposable injective is isomorphic to
$E(\cO_Z)$ for some closed reduced and irreducible subscheme 
$Z$ of $X$ \cite[Th\'eor\`eme 1, p. 443]{G}.
It therefore suffices to show that every
coherent submodule of each $E(\cO_Z)$ is good.

Fix a closed reduced and irreducible subscheme $Z \subset X$, and a
coherent $\cO_X$-submodule $\cM \subset E(\cO_Z)$. 
Let $z$ denote the generic point of $Z$, and let $\cO_z$ denote 
the stalk of $\cO_X$ at $z$.
There is a morphism $f:\Spec \cO_{z} \to X$
with the following
properties: the inverse image functor $f^*$ is exact, the direct image
functor $f_*$ is fully faithful and exact, and has a right adjoint 
$f^!$. Because $f_*$ is fully faithful the counit $f^*f_* \to
\id_{\Spec \cO_{z}}$ is an isomorphism.

Let $i:Z \to X$ be the inclusion. Let $\cE$ be the constant sheaf
on $Z$ having sections the function field of $Z$. Then $i_*\cE$ is
an essential extension of $\cO_Z$, so $E(\cO_Z)=E(i_*\cE)$.
But $i_*\cE$ is also gotten by applying $f_*$ to the residue field
of $\cO_{z}$, so the unit $i_*\cE \to f_*f^*(i_*\cE)$ is an 
isomorphism. 
However, $f_*$ sends injectives to injectives because it is 
right adjoint to an exact functor, so
if $\cF$ is an injective envelope of $f^*(i_*\cE)$ in $\Mod
\cO_{z}$, $f_*\cF$ is an injective quasi-coherent $\cO_X$-module 
containing a copy of $i_*\cE$. 
Thus $E(\cO_Z) \cong f_*\cF$. There is a 
surjective map $\cO_{z}^{(I)} \to \cF$ from a suitably large direct
sum of copies of $\cO_{z}$, and therefore an epimorphism 
$f_*(\cO_{z}^{(I)}) \to f_*\cF$. Since $f_*$ has a right adjoint it
commutes with direct sums, so we obtain an epimorphism
$(f_*\cO_{z})^{(I)} \to f_*\cF$. Because $\Qcoh X$ is locally
noetherian, every coherent $\cO_X$-submodule of $f_*\cF$ is
therefore an epimorphic image of a coherent submodule of
$f_*\cO_{z}^{(I)}$. However, $f_*\cO_{z}$ is an $\cO_X$-submodule of
$\cK$, so every coherent $\cO_X$-submodule of it is good.
It follows that every coherent submodule of $f_*\cF$ is good.
Hence $\cM$ is good.
\end{pf}

\begin{corollary}
\label{cor.integral.integral}
Let $X$ be a noetherian scheme. 
Then $X$ is integral in the usual sense if and only if it is
integral in the sense of Definition \ref{defn.integral}.
In that case, $\cE_X$ is
isomorphic to the constant sheaf $\cK$ with sections the function
field of $X$.
\end{corollary}
\begin{pf}
Let $X$ be integral in the usual sense of algebraic geometry.
By \cite[Chapitre VI]{G}, $\cK$ is an indecomposable injective.
It is also clear from Gabriel's classification of the indecomposable
injectives in $\Qcoh X$ that $\cK$ is the only indecomposable
injective of Krull dimension equal to $\dim X$. It therefore
follows from Proposition \ref{prop.intgl.scheme}, that $X$ is integral 
in our sense. Furthermore,
$\cE_X=\cK$, and the endomorphism ring of $\cK$ is $k(X)$, so 
function field and generic point in our sense agree with the usual
notions.

Conversely, suppose that $X$ is integral in the sense of
Definition \ref{defn.integral}.
By \cite{G}, $\cE_X \cong E(\cO_Z)$ for some closed reduced and
irreducible subscheme $Z$ of $X$. In particular, $\rank \cO_Z=1$.

We will show that every coherent $\cO_X$-submodule of $\cE_X$ is an
$\cO_Z$-module. It will then follow that the same is true of every
coherent subquotient of a finite direct sum of copies of $\cE_X$. In
particular, $\cO_X$ will be an $\cO_Z$-module, whence $Z=X$, and the
proof is complete.

It suffices to prove that every coherent submodule of $E(\cO_Z)$ 
containing $\cO_Z$ is an $\cO_Z$-module. Let $M$ be such a
submodule. If $W$ denotes the support of $M/\cO_Z$, then
$M/\cO_Z$ is annihilated by some
power of $\cI_W$, the ideal cutting out $W$. Hence
$M\cI_W^n\cI_Z=0$ for $n \gg 0$. If $M \cI_Z=0$, then $M$ is an
$\cO_Z$-module, so we may suppose that $M\cI_Z$ is non-zero. 
Hence $M\cI_Z$ 
has non-zero intersection with the essential submodule $\cO_Z$ of
$E(\cO_Z)$, so $\cI^n_W$ annihilates a non-zero ideal of
$\cO_Z$. But $Z$ is integral, so this can only happen if $\cI_W^n
\subset \cI_Z$; hence $Z \subset W$.

On the other hand the injective envelope of $M/\cO_Z$ is a direct
sum of indecomposable injectives, so a direct sum of copies of
$E(\cO_{W_i})$ for various closed integral subschemes $W_i$ of $X$.
Since $Z$ is contained in the support of $M/\cO_Z$, and every
non-zero coherent submodule of $E(\cO_{W_i})$ has support equal to
$W_i$, $Z$ is contained in the union of the $W_i$s. Since $Z$ is
integral it must be contained in one of the $W_i$s. Hence $\cO_Z$
is a quotient of $\cO_{W_i}$ for some $i$, and we deduce that
$\rank \cO_{W_i} \ge 1$. It follows that the rank of
$E(M/\cO_Z)=1$.
Hence by Proposition \ref{prop.essl.extn}, $\rank M/\cO_Z =1$. This
contradicts the fact that $\rank \cE_X/\cO_Z=0$, so we conclude
that $M\cI_Z=0$.
Hence $M$ is an $\cO_Z$-module, as required.
\end{pf}

\begin{proposition}
\label{prop.goldie}
Let $R$ be a right noetherian ring and let $X$ be the affine 
space with coordinate ring $R$. 
If $R$ is prime,  then $X$ is integral.
\end{proposition}
\begin{pf}
By Goldie's theorem, the ring of fractions of $R$ is a matrix 
ring over a division ring, say $D$. Furthermore, that matrix ring is an
injective envelope of $R$ as a right $R$-module.
Let $\cE$ be a simple right ideal of that matrix ring.
The endomorphism ring of $\cE$ as an $R$-module is the same as its
endomorphism ring as a module over the matrix ring, so is equal to
the division ring $D$.
Since $R$ embeds in a finite direct sum of copies of $\cE$, and
is a generator in $\Mod R$, every noetherian right
$R$-module is a subquotient of a direct sum of copies of $\cE$.
\end{pf}

In Proposition \ref{prop.goldie}, 
the function field of $X$ is the division ring $D$ that appears in
Goldie's Theorem.

It is not the case that a right noetherian ring $R$ is prime if
and only if $\Mod R$ is integral.
For example, the ring of upper triangular matrices over a field is
integral in our sense. 
However, it is easy to see that if $X$ is affine and integral, then
its coordinate ring is prime if and only if $\cE_X$ is a prime
$X$-module in the sense of \cite[Definition 4.3]{Sm}.
Proposition \ref{prop.goldie.2} also gives a criterion which implies 
that the coordinate ring of an integral affine space is prime.

Our notion of integral is not an invariant of the derived category.
For example, let $A$ be the path algebra of the
quiver $\bullet \to \bullet \to \bullet$
and $A'$ the path algebra of the 
 quiver $\bullet \leftarrow \bullet \rightarrow \bullet$. The
derived categories of modules over $A$ and $A'$ are equivalent.
By listing the three indecomposable injectives over each algebra it
is clear that $\Mod A$ is integral, but $\Mod A'$ is not.
In fact, the path algebra of a quiver without loops is integral if
and only if it has a unique sink.
We are grateful to D. Happel for these observations.

\medskip

We now show that a non-commutative analogue of a projective scheme
is integral if it has a homogeneous coordinate ring that is prime
and noetherian.

\begin{definition}
[Verevkin \cite{Ver}, Artin and Zhang \cite{AZ}]
\label{defn.proj.space}
Let $A$ be an $\NN$-graded $k$-algebra such that $\dim_k A_n <
\infty$ for all $n$. Define $\GrMod A$ to be the category
of $\ZZ$-graded $A$-modules with morphisms the $A$-module homomorphisms
of degree zero. We write $\Fdim A$ for the full subcategory of direct limits 
of finite dimensional modules. We define the quotient category
$$
\Tails A=\GrMod A/\Fdim A,
$$
and denote by $\pi$ and $\omega$ the quotient functor and its right adjoint.
The {\sf projective space} $X$ with {\sf homogeneous coordinate ring}
$A$ is defined by $\Mod X:=\Tails A$. 
\end{definition}

\begin{theorem}
\label{thm.intgl.tails}
Let $A$ be prime noetherian locally finite $\NN$-graded $k$-algebra. 
Suppose that $\dim_k A= \infty$. 
Suppose further that the graded ring of fractions $\Fract_{gr} A$
contains an isomorphic copy of $A(n)$ for every integer $n$.
Then the projective space with homogeneous coordinate ring $A$
is locally noetherian and integral. Its function field
is the degree zero component of $\Fract_{gr} A$.
\end{theorem}
\begin{pf}
Define $X$ by $\Mod X=\Tails A$.
Since $\Mod X$ is a quotient of a locally noetherian
category it is locally noetherian.

It is well-known that the injective envelope of $A$ in $\GrMod A$
is its graded ring of fractions, say $E=\Fract_{gr} A$.
Let $\cE=\pi E$ be its image in $\Mod X$.
Since $A$ is prime and has infinite dimension, zero is the only
finite dimensional graded submodule of it. The same is true of $E$,
so $\cE$ is injective in $\Mod X$.

To show that $X$ is integral it only remains to
show that every noetherian $X$-module is a subquotient of a
finite direct sum of copies of $\cE$. If $\cM$ is a noetherian
$X$-module, then $\cM \cong \pi M$ for some noetherian $A$-module
$M$. Now
$M$ is a quotient of a finite direct sum of shifts $A(n)$ for
various integers $n$, so $\cM$ is a quotient of a 
finite direct sum of various twists $\cO_X(n)=\pi A(n)$. However, each
$A(n)$ embeds in $E$, so each $\cO_X(n)$ embeds in $\cE$. Thus
$\cM$ is a subquotient of a finite direct sum of copies of $\cE$.

Finally, 
$$
k(X)=\Hom_X(\cE,\cE)=\Hom_X(\pi E,\pi E) \cong \Hom_{\Gr
A}(E,\omega\pi E).
$$
However, since zero is the only finite dimensional submodule of $E$
and $E$ is injective, $\omega\pi E \cong E$. Hence
$$
k(X) \cong \Hom_{\Gr A}(E,E) \cong (\Fract_{gr} A)_0,
$$
as claimed.
\end{pf}

The hypothesis in Theorem \ref{thm.intgl.tails} that
$\Fract_{gr}A$ contain a copy of each $A(n)$ is necessary because if
$A=k[x^2]$ 
with $\deg x=1$, then $X \cong \Spec k^2$.
This hypothesis holds if $A$ has a
regular element in all sufficiently high degrees. In particular, if
$A$ is a domain generated in degree one, then $X$ is integral. 
Thus, the quantum planes of
Artin-Tate-Van den Bergh are integral, as are the other standard
non-commutative analogues of the projective spaces $\PP^n$.

Van den Bergh has defined the notion of the blowup at a closed
point on a non-commutative surface \cite{vdB}. The exceptional fiber,
$E$ say, is sometimes, but not always, a projective line. 
Nevertheless it is always integral. For example, when
$\Mod E=\GrMod k[x]$ its big injective is $k[x,x^{-1}]$ and its
function field is $k$. In the other
cases $\Mod E$ is of the form $\Tails k[x,y]$ where $k[x,y]$ 
is the commutative polynomial ring with $\deg x=1$ and 
$\deg y=n<\infty$, and its integrality is guaranteed by Theorem
\ref{thm.intgl.tails}. In these cases the function field of $E$
is the rational function field $k(y/x^n)$.


\section{Properties of integral spaces}

An integral scheme has several properties that we
might expect a non-commutative integral space to have. 
For example, every non-empty open subscheme of a noetherian 
integral scheme is dense because it 
contains the generic point. To get a non-commutative version of
this we must first introduce analogues of ``open subspace'' and
``closure''. This is done in \cite{Sm}, but we recall the
definitions here.

\begin{definition}
\label{defn.open}
Let $X$ be a non-commutative space. 
A {\sf weakly open subspace}, say $U$,
of $X$ is a full subcategory $\Mod U$ of $\Mod X$ such that the
inclusion functor $\a_*:\Mod U \to \Mod X$ has an exact left adjoint
$\a^*$.
\end{definition}

For example, the generic point of an integral space is a weakly open
subspace.

\begin{definition}
\label{defn.weakly.closed}
A {\sf weakly closed subspace} $W$ of a non-commutative space
$X$ is a full subcategory $\Mod W$ of $\Mod X$ 
that is closed under subquotients and isomorphisms, and for which 
the inclusion functor $\a_*:\Mod W \to \Mod X$ has 
a right adjoint.
We write $\a:W \to X$ for the weak map corresponding to $\a_*$. 
\end{definition}

Let $\a:W \to X$ be the inclusion of a weakly closed subspace.
Then $\Mod W$ is a Grothendieck category and is
locally noetherian if $\Mod X$ is. Because $\Mod W$ is closed 
under subquotients, $\a_*$ is an exact functor.
Because $\a_*$ has a right adjoint it commutes with direct sums.
Further information about weakly closed subspaces 
can be found in \cite{Sm}. 

The requirement in the definition of an integral space that every
$X$-module be a subquotient of a direct sum of copies of $\cE_X$ is
equivalent to the requirement that $X$ is the only weakly closed
subspace having $\cE_X$ as a module over it.

\medskip

Let $U$ and $Z$ be respectively a weakly open and a weakly closed
subspace of $X$. We say that $Z$ {\sf contains} $U$ if $\Mod U$ is
contained in $\Mod Z$. In other words, if $\a:U \to X$ and $\d:Z \to
X$ are the inclusions, then $U$ is contained in $Z$ if and only if
there is a weak map $\ve:U \to Z$ such that $\d\ve=\a$. In this case,
$U$ becomes a weakly open subspace of $Z$ because $\a^*\d_*$ 
is an exact left adjoint to $\ve_*$: if $M \in \Mod Z$ and $N \in
\Mod U$, then
$$
\Hom_Z(M,\ve_*N) = \Hom_X(\d_*M,\d_*\ve_*N)=\Hom_X(\d_*M,\a_*N) 
\cong \Hom_U(\a^*\d_*M,N).
$$

\begin{definition}
\cite{Sm} 
If $U$ is a weakly open subspace of a locally noetherian
space $X$ its {\sf weak closure}, denoted $\Uol$,
is the smallest weakly closed subspace of $X$ that contains $U$. 
\end{definition}

This makes sense because the intersection of two
weakly closed subspaces is a weakly closed subspace.
If $\a:U \to X$ is the inclusion, then
$\Mod \Uol$ consists of all subquotients of $X$-modules
of the form $\a_*N$ as $N$ ranges over $\Mod U$.
More details about weak closure can be found in \cite{Sm}.


\begin{lemma}
\label{lem.gen.pt.dense}
If $\eta$ is the generic point of an integral space $X$, then
$\overline{\eta}=X$.
\end{lemma}
\begin{pf}
If $Z$ is a weakly closed subspace of $X$ containing $\eta$, then
$\cE_X$ belongs to $\Mod Z$. Since $\Mod Z$ is closed under
subquotients and direct sums, every $X$-module belongs to $\Mod Z$,
showing that $Z=X$.
\end{pf}

\begin{lemma}
\label{lem.dense.pt}
Let $p$ be a weakly open point in a locally noetherian space $X$. 
That is, $p$ is a weakly open subspace of $X$
and $\Mod p=\Mod D$ for some division ring $D$. If $\pol=X$, then
$X$ is integral, $p$ is its generic point, and $k(X)=D$.
\end{lemma}
\begin{pf}
Let $\a:p \to X$ denote the inclusion. The big injective in $\Mod
p$ is $D$. Since $\a_*$ is right adjoint to an exact functor,
$\cE:=\a_*D$ is an injective $X$-module. Using the adjoint pair 
$(\a^*,\a_*)$ it is easy to see that $\cE$
is indecomposable because $D$ is, and its endomorphism ring is the
same as that of $D$, namely $D$. Furthermore, if $M$ is an $X$-module,
it is a subquotient of $\a_*N$ for some $p$-module $N$ because
$\pol=X$. But $N$ is a direct sum of copies of $D$, and $\a_*$ 
commutes with direct sums \cite[Cor. 1, p. 379]{G}, so $M$ is
a subquotient of a direct sum of copies of $\cE$. Hence $X$ is
integral.

To see that $p$ is the generic point of $X$ it suffices to show
that $\a^*$ vanishes on the torsion modules. However, if $M$ is
torsion, then $0=\Hom_X(M,\cE) \cong \Hom_p(\a^*M,D)$, whence
$\a^*M=0$.
\end{pf}

\begin{proposition}
\label{prop.dense.intgl}
Let $U$ be a weakly open subspace of a locally noetherian space $X$. 
Suppose that $U$ is integral and $\Uol=X$. If the inclusion $U
\to X$ is an affine map, then $X$ is integral and $k(X)=k(U)$.
\end{proposition}
\begin{pf}
The notion of an affine map is defined in \cite{Sm}; the important
point here is that if $\a:U \to X$ denotes the inclusion, then
$\a_*$ is exact. Let $\cE_U$ be the big injective in $\Mod U$. Since
$\a_*$ is right adjoint to an exact functor, $\a_*\cE_U$ is an
injective $X$-module. It is also indecomposable, and its
endomorphism ring is equal to $\End_U \cE_U$. 

It remains to show that every $X$-module is a subquotient of
a direct sum of copies of $\a_*\cE_U$. 
Let $P \in \Mod U$.
Since $U$ is integral, $P \cong B/A$ for some $U$-submodules $A
\subset B \subset \cE_U^{(I)}$ and some index set $I$. Since $\a_*$
is exact, $\a_*P \cong (\a_*B)/(\a_*A)$; since $\a_*$ commutes with
direct sums we have $X$-submodules $\a_*A \subset \a_*B \subset
(\a_*\cE_U)^{(I)}$; thus $\a_*P$ is a subquotient of a direct sum
of copies of $\a_*\cE_U$. But $\Uol=X$, so every
$X$-module is a subquotient of $\a_*P$ for some $P \in \Mod U$. 
The result now follows.
%
\end{pf}

Proposition \ref{prop.dense.intgl} applies to the situation where
one has an affine space and embeds it in a projective space by
adding an effective divisor at 
infinity (see \cite[Section 8]{Sm})---if
the affine space is integral, so is the projective space, and their
function fields coincide.

\medskip
Let $W$ be a weakly closed subspace of a locally noetherian space
$X$. Its complement $X \backslash W$ is defined 
in \cite[Section 6]{Sm}. In particular, $X \backslash W$ is a
weakly open subspace of $X$, and every weakly open subspace arises as
such a complement.

%
%
%
%
%
%

\begin{proposition}
\label{prop.good.subspace}
Let $X$ be an integral space and $W$ a weakly closed subspace.
Suppose that $W \ne X$.
If $\cE_X$ does not contain a non-zero $W$-submodule, then
\begin{enumerate}
\item{}
$\eta$ belongs to $X \backslash W$ and $\overline{X \backslash W}=X$;
\item{}
$X \backslash W$ is integral and $k(X \backslash W)=k(X)$.
\end{enumerate}
\end{proposition}
\begin{pf}
Let $\a:X \backslash W \to X$ denote the inclusion.
Let $\tau:\Mod X \to \Mod X$ denote the functor that 
is the kernel of the natural transformation $\id_X \to \a_*\a^*$.
There is an exact sequence
$$
0 \to \tau \cE_X \to \cE_X \to \a_*\a^*\cE_X \to R^1\tau \cE_X \to 0.
$$
By hypothesis, $\tau \cE_X=0$.  Since $\cE_X$ is injective, 
$R^1\tau \cE_X=0$. Hence $\cE_X \cong \a_*\a^*\cE_X$. 
It follows that the generic point of $X$ belongs to $X \backslash W$.
More formally, if $j:\eta \to X$ is the inclusion, 
then there is a map $\c:\eta \to X \backslash W$ 
such that $j=\a\c$ ( this is straightforward, though it can also 
be seen as a special case of \cite[Proposition 6.1]{Sm}).
By Lemma \ref{lem.gen.pt.dense}, the weak closure
of $X \backslash W$ is $X$. This proves (1).

Because $\tau\cE_X=0$, $\a^*\cE_X$ is an injective 
$X \backslash W$-module. It is an indecomposable injective because
$$
\Hom_{X \backslash W}(\a^*\cE_X,\a^*\cE_X) \cong
\Hom_X(\cE_X,\a_*\a^*\cE_X) = \Hom_X(\cE_X,\cE_X)
$$
is a division ring.
If $\cM$ is a noetherian $X \backslash W$-module, 
then $\cM=\a^*M$ for some 
noetherian $X$-module $M$. There is a noetherian submodule $L$ of 
$\cE_X^{\oplus n}$ and an epimorphism $L \to M$. Hence $\a^*L$ is a
noetherian submodule of $\a^*\cE_X^{\oplus n}$ and there is an 
epimorphism $\a^*L \to \a^*M$. Thus $X \backslash W$ is integral.
\end{pf}


%

We define the {\sf empty space} $\phi$ by declaring $\Mod \phi$ to
be the zero category; that is, the abelian category having only 
one object and one morphism.
Part (1) of Proposition \ref{prop.good.subspace} can now be
rephrased as follows.
If $W_1$ and $W_2$ are non-empty weakly closed subspaces of an
integral space $X$ such that $\cE_X$ contains  neither a
non-zero $W_1$-module nor a non-zero $W_2$-module, then 
$(X \backslash W_1) \cap (X \backslash W_2) \ne \phi$.
By \cite[Section 6]{Sm}, this intersection is equal to $X\backslash
(W_1 \cup W_2)$, so we deduce that $W_1 \cup W_2 \ne X$.

\section{Dimension Functions}

M. Van den Bergh has suggested that a dimension function should 
play a prominent role in non-commutative geometry. 

In an earlier version of this paper our definition of integrality 
required the big injective to be critical with respect to a 
dimension function. We are grateful to the referee for suggesting
that this was unnecessary. Nevertheless, since dimension functions
play an important role in non-commutative algebra and geometry it
is useful to examine the connection. 

\begin{definition}
\label{defn.dim}
Let $X$ be a locally noetherian space.
A {\sf dimension function} on $X$ is a function $\d:\Mod X \to 
\RR_{\ge 0} \cup \{-\infty,\infty\}$ satisfying the following conditions:
\begin{itemize}
\item{}
$\d(0)=-\infty$;
\item{}
if $0 \to L \to M \to N \to 0$ is exact, then
$\d(M)=\max\{\d(L),\d(N)\}$;
\item{}
$\d(M)=\max\{\d(N) \; | \; N \hbox{ is a noetherian submodule of
$M$}\}$;
\item{}
if $\s$ is an auto-equivalence of $\Mod X$, then $\d(M^\s)=\d(M)$.
\end{itemize}
We define the dimension of $X$, $\dim X$, to be the maximum of
$\d(M)$ as $M$ ranges over all $X$-modules.
\end{definition}

{\bf Remarks. 1.}
We will make no use in this paper of the condition that $\d$
is invariant under auto-equivalences. 

{\bf 2.}
A dimension function $\d$ determines various localizing subcategories
of $X$. If $d \in \RR_{\ge 0} \cup \{\infty\}$, we 
write $\Mod_{\le d} X$
and $\Mod_{< d} X$ for the full subcategories of $\Mod X$ 
consisting of those $M$ such that $\d(M) \le d$ and $\d(M)<d$
respectively.  These are localizing subcategories because
$\d(\sum_j N_j)=\max_j\d(N_j)$.
One can specify the dimension function simply by
specifying these localizing subcategories.

{\bf 3.}
The notion of Krull dimension as defined by Gabriel in
\cite{G} is a dimension function. It is defined inductively:
$\Mod_{< 0} X$ consists of only the zero module, and
for each integer $n \ge 0$,
$\Mod_{\le n} X / \Mod_{<n} X$ consists of all direct 
limits of artinian modules in $\Mod X/\Mod_{< n} X$.

The version of Krull dimension defined
using posets that appears in \cite[Chapter 6]{MR}, does {\em not}
satisfy our definition of dimension function. In fact, it is not even
defined for all modules, and does not lead to an ascending chain of
localizing subcategories.
Thus, we always use Gabriel's version of Krull dimension.

{\bf 4.}
If $X$ is a noetherian scheme, then the Krull dimension of a
coherent $\cO_X$-module is equal to the dimension of its support.

{\bf 5.}
Each of the localizing subcategories described above determines a 
subgroup of $K_0(X)$, and in this way one obtains a filtration of
$K_0(X)$.

{\bf 6.}
If $X$ is a locally noetherian space with a dimension function $\d$, 
then every weakly closed subspace of $X$ is locally noetherian, and
it inherits the dimension function. The dimension of such a 
subspace is the maximum of the dimensions of its noetherian modules.

\begin{definition}
An $X$-module $M$ is {\sf $d$-critical} if $\d(M)=d$ and
$\d(M/N)<d$ for all non-zero submodules $N$ contained in $M$.
We say that $M$ is {\sf $d$-pure} if $\d(N)=d$ for 
all its non-zero submodules $N$.
The {\sf $d$-length} of an $X$-module $M$ is its length in 
$\Mod X/\Mod_{< d} X$. It is denoted by $\ell_d(M)$, and it may
take the value $\infty$.
\end{definition}

Let $X$ be a noetherian scheme with Krull dimension as 
the dimension function.
If $Z$ is a closed subscheme of $X$, then $\cO_Z$ is critical in 
$\Qcoh X$ if and only if $Z$ is an integral subscheme of $X$.

The function $\ell_d(-)$ is additive on short exact sequences.
One sees this by passing to the quotient
category $\Mod X/\Mod_{< d} X$ and using the fact that
the usual notion of length is additive.
Because $\ell_d$ is additive, a $d$-critical module is uniform
(i.e., two non-zero submodules of it have non-zero
intersection). Hence an injective envelope of a $d$-critical module
is indecomposable.

If $M$ is a noetherian module of dimension $d$, then $M$ has an
$d$-critical quotient module, namely $M/N$, where $N$ is a submodule
of $M$ maximal subject to the condition that $\d( M/N)=d$.

A $d$-critical module is $d$-pure.
A $d$-pure module is critical if and only if its 
$d$-length is one.

\begin{proposition}
Let $X$ be a locally noetherian space.
Suppose that $\cE$ is an indecomposable injective such that
every $X$-module is a subquotient of a direct sum of copies of
$\cE$. If $\cE$ is $d$-critical with respect some dimension
function, then 
\begin{enumerate}
\item
$X$ is integral of dimension $d$ and $\cE$ is the big injective;
\item
$M$ is torsion if and only if $\d(M)<d$;
\item
$\ell_d(M)=\rank M$.
\end{enumerate}
\end{proposition}
\begin{pf}
(1)
If $M$ is a non-zero submodule of $\cE$, then
$\d(\cE/M)<\d(\cE)$, whence $\Hom_X(\cE/M,\cE)=0$.
It follows that the endomorphism ring of $\cE$ is a division ring. 
Hence $X$ is integral.
Since an $X$-module is a subquotient of a direct sum of copies of
$\cE$ its dimension is at most $d$. Hence $\dim X=d$.

(2)
If $\d(M)<d$, then $\Hom_X(M,\cE)=0$ because $\cE$ is $d$-critical,
and $M$ is torsion. To prove the converse it suffices to show
if $M$ is a noetherian module such that $\d(M)=d$, then $M$ is not
torsion. Suppose to the contrary that there is such an $M$ which is
torsion. Then $M$ has a $d$-critical quotient $\Mol$. This is also
torsion, and so is its injective envelope $E(\Mol)$ by Proposition
\ref{prop.torsion}. By Proposition \ref{prop.red+irr}, there is a
non-zero map $\varphi:\cE \to E(\Mol)$. Since $E(\Mol)$ is torsion,
$\varphi$ is not monic. Since $\cE$ is $d$-critical, $\d(\im
\varphi)<d$. Hence $\d(\im \varphi \cap \Mol)<d$. But $\im \varphi
\cap \Mol \ne 0$, so this contradicts the fact that $\Mol$ is
$d$-critical. We conclude that $M$ can not be torsion.

(3)
By (2), $\Mod_{<d}X$ consists of the torsion
modules, whence $\Mod X/\Mod_{<d}X=\Mod \eta$, where $\eta$ is the
generic point of $X$. 
The remark after Proposition \ref{prop.gen.pt}
implies that $\rank M=\ell_d(M)$.
\end{pf}

\begin{proposition}
\label{prop.goldie.2}
Let $X$ be an integral locally noetherian affine space 
with coordinate ring $R$. 
Suppose there is a dimension function $\d$ such that 
$\d(M \otimes_R I) \le
\d(M)$ for all noetherian modules $M$ and all two-sided ideals $I$.
If $\cE_X$ is critical with respect to $\d$,
then $R$ is prime.
\end{proposition}
\begin{pf}
Since $X$ is locally noetherian, $R$ is right noetherian.
By \cite[Prop. 3.9]{JTS}, the condition on $\d$ ensures that the
annihilator of a critical right $R$-module is a prime ideal. In
particular, $\Ann \cE_X$ is a prime ideal. But $R$ itself is a
subquotient of a finite direct sum of copies of $\cE_X$, so the
annihilator of $\cE_X$ is zero. Hence $R$ is prime.
\end{pf}

We expect there is a dimension function for right noetherian rings 
satisfying the hypothesis in Proposition \ref{prop.goldie.2}.
For many two-sided noetherian rings, such as factors of enveloping
algebras,  Gelfand-Kirillov dimension satisfies the hypothesis.

\smallskip

Every proper closed subscheme of an integral noetherian scheme $X$ 
has strictly smaller dimension than $X$.
For non-commutative spaces Krull dimension does not necessarily
have this property---for example, take the ring of upper 
triangular matrices over a field.


We now pick out a better behaved class of weakly closed subspaces.

\begin{definition}
\label{defn.good}
Let $\d$ be a dimension function on $X$.
A weakly closed subspace $W$ of $X$
is {\sf good} if whenever $0 \to L \to M \to N \to 0$ is an 
essential extension of a $W$-module $L$ by an $X$-module $N$ such that
$\d(N)<\d(L)$, then $M \in \Mod W$.
\end{definition}

A subspace can be good with respect to one dimension function but
not good with respect to another.

If $X$ is integral and $W \subset X$ is a proper weakly closed
subspace, then $\dim W < \dim X$ if $W$ is good. Hence we have the
following result.

\begin{lemma}
\label{lem.maxl.chains}
Let $X$ be an integral space, and suppose that $\d(M)\in \NN$ 
for all $M \ne 0$. If 
$$
\phi \ne W_0 \subset W_1 \subset \ldots \subset W_d
$$
is a chain of distinct good integral subspaces of $X$, then $d \le \dim
X$.
\end{lemma}

\begin{example}
\label{eg.triang}
Let $R$ be the ring of lower triangular $2\times 2$ matrices
over a field. Let $\cO_p$ and $\cO_q$ be the two simple 
right $R$-modules with $\cO_p$ the projective one.
There are closed points, $p$ and $q$, defined by $\Mod p$ consists of 
all direct sums of copies of $\cO_p$; $\Mod q$ is defined similarly 
(closed points are defined in \cite{Sm}).
There is a non-split exact sequence
$0 \to \cO_p \to \fp \to \cO_q \to 0,$
where $\fp$ is the annihilator of $\cO_p$. 
The indecomposable injectives are $\cO_q$ and $\fp$. 

Since $\End_R(\fp) \cong k$ and 
$R \cong \fp \oplus \fp/\cO_p$, every $R$-module is a
subquotient of a direct sum of copies of $\fp$. Therefore
$X$ is integral, $\fp$ is the big injective, 
and the function field of $X$ is $k$.
If $j:\eta \to X$ is the inclusion of the generic point, then
$j_*(\Mod \eta)$ consists of all direct sums of copies of $\fp$. 
We also note that $\eta=X \backslash q$.

There are several ways in which $X$ does not behave like an
integral scheme.
The inclusion $X \backslash p \to X$ sends $X \backslash p$
isomorphically onto $q$, so $X \backslash p$ is both open and
closed in $X$. 
In particular, $\overline{X \backslash p} \ne X$.
Furthermore, if we view $\eta$ as an open subspace of $X$, then 
$\eta \cap (X \backslash p) = \phi$.
Finally $p$ is a proper closed subspace of $X$ having the
same Krull dimension as $X$.
\end{example}

\end{document}